\documentclass[12pt, twoside, leqno]{article}

\usepackage{amsmath,amsthm}
\usepackage{amssymb}

\usepackage{enumerate}

\newtheorem{thm}{Theorem}[section]
\newtheorem{cor}[thm]{Corollary}
\newtheorem{lem}[thm]{Lemma}

\theoremstyle{definition}

\numberwithin{equation}{section}

\begin{document}


\baselineskip=17pt


\title{Improved lower bound for the number of unimodular zeros of self-reciprocal polynomials with coefficients in a finite set}

\author{Tam\'as Erd\'elyi\\
Department of Mathematics\\
Texas A\&M University\\
College Station, Texas 77843}

\date{January 23, 2019}

\maketitle


\renewcommand{\thefootnote}{}



\renewcommand{\thefootnote}{\arabic{footnote}}
\setcounter{footnote}{0}


\begin{abstract}
Let $n_1 < n_2 < \cdots < n_N$ be non-negative integers. In a private communication Brian Conrey asked
how fast the number of real zeros of the trigonometric polynomials $T_N(\theta) = \sum_{j=1}^N {\cos (n_j\theta)}$
tends to $\infty$ as a function of $N$. Conrey's question in general does not appear to be easy.
Let ${\mathcal P}_n(S)$ be the set of all algebraic polynomials of degree at most $n$ with each of their
coefficients in $S$. For a finite set $S \subset {\mathbb C}$ let $M = M(S) := \max\{|z|: z \in S\}$.
It has been shown recently that if $S \subset {\mathbb R}$ is a finite set and $(P_n)$ is a
sequence of self-reciprocal polynomials $P_n \in {\mathcal P}_n(S)$ with $|P_n(1)|$ tending to $\infty$,
then the number of zeros of $P_n$ on the unit circle also tends to $\infty$. In this paper
we show that if $S \subset {\mathbb Z}$  is a finite set, then every self-reciprocal polynomial
$P \in {\mathcal P}_n(S)$ has at least $$c(\log\log\log|P(1)|)^{1-\varepsilon}-1$$
zeros on the unit circle of ${\mathbb C}$ with a constant $c > 0$ depending only on $\varepsilon > 0$ and
$M = M(S)$.
Our new result improves the exponent $1/2 - \varepsilon$ in a recent result by Julian Sahasrabudhe to
$1 - \varepsilon$. Sahasrabudhe's new idea \cite{Sah} is combined with the approach used in \cite{Erd8} 
offering an essencially simplified way to achieve our improvement. We note that in both Sahasrabudhe's 
paper and our paper the assumption that the finite set $S$ contains only integers is deeply exploited.
\end{abstract}

\newpage
\section{Introduction and Notation.}

Research on the distribution of the zeros of algebraic polynomials has a
long and rich history. In fact, most of the papers \cite{AndBla}--\cite{TotVar} in our list of references
are just some of the papers devoted to this topic. The study of the number of
real zeros of trigonometric polynomials and the number of unimodular zeros
(that is, zeros lying on the unit circle of the complex plane) of algebraic
polynomials with various constraints on their coefficients is the subject of
quite a few of these. We do not try to survey these in our introduction.


Let $S \subset {\mathbb C}$. Let ${\mathcal P}_n(S)$ be the set of all algebraic polynomials of
degree at most $n$ with each of their coefficients in $S$. A polynomial $P$ of the form
\begin{equation}
P(z) = \sum_{j=0}^n{a_jz^j}\,, \qquad a_j \in {\mathbb C}\,, 
\end{equation}
is called conjugate-reciprocal if
\begin{equation}
\overline{a}_j = a_{n-j}\,, \qquad j=0,1,\ldots,n\,.
\end{equation}
A polynomial $P$ of the form (1.1) is called plain-reciprocal or self-reciprocal if
\begin{equation}
a_j = a_{n-j}, \qquad j=0,1,\ldots,n\,. 
\end{equation}
If a conjugate reciprocal polynomial $P$ has only real coefficients, then it is obviously
plain-reciprocal. We note also that if
$$P(z) = \sum_{j=0}^{2n}{a_j z^j}\,, \qquad a_j \in {\mathbb C}\,,$$
is conjugate-reciprocal, then there are $\theta_j \in {\mathbb R}$, $j=1,2,\ldots\,n$, such that
$$T(t) := P(e^{it})e^{-int} = a_n + \sum_{j=1}^n{2|a_{j+n}|\cos(jt + \theta_j)}\,.$$
If the polynomial $P$ above is plain-reciprocal, then
$$T(t) := P(e^{it})e^{-int} = a_n + \sum_{j=1}^n{2a_{j+n}\cos(jt)}\,.$$
In this paper, whenever we write ``$P \in {\mathcal P}_n(S)$ is conjugate-reciprocal"
we mean that $P$ is of the form (1.1) with each $a_j \in S$ satisfying (1.2).
Similarly, whenever we write ``$P \in {\mathcal P}_n(S)$ is self-reciprocal" we mean that
$P$ is of the form (1.1) with each $a_j \in S$ satisfying (1.3). This is going to be
our understanding even if the degree of $P \in {\mathcal P}_n(S)$ is less than $n$.
It is easy to see that $P \in {\mathcal P}_n(S)$ is self-reciprocal and $n$ is odd, then
$P(-1)=0$.
Associated with an algebraic polynomial $P$ of the form (1.1) we introduce the numbers
$$\text{\rm NC}(P) := |\{j \in \{0,1,\ldots,n\}: a_j \neq 0\}|\,.$$
Here, and in what follows $|A|$ denotes the number of elements of a finite set $A$.
Let $\text{\rm NZ}(P)$ denote the number of real zeros (by counting multiplicities) of an algebraic
polynomial $P$ on the unit circle. Associated with a trigonometric polynomial
$$T(t) = \sum_{j=0}^n{a_j \cos(jt)}$$
we introduce the numbers
$$\text{\rm NC}(T) := |\{j \in \{0,1,\ldots,n\}: a_j \neq 0\}|\,.$$
Let $\text{\rm NZ}(T)$ denote the number of real zeros (by counting multiplicities) of a real 
trigonometric polynomial $T$ in a period (of length $2\pi$). Let $\text{\rm NZ}^*(T)$ denote the 
number of sign changes of a real trigonometric polynomial $T$ in a period (of length $2\pi$).
The quotation below is from \cite{BorErdFerLoc}.
``Let $0 \leq n_1 < n_2 < \cdots < n_N$ be integers. A cosine polynomial of the form
$T(\theta) = \sum_{j=1}^N {\cos (n_j\theta)}$ must have at least one real zero in
a period. This is obvious if $n_1 \neq 0$, since
then the integral of the sum on a period is $0$.
The above statement is less obvious if $n_1 = 0$, but for
sufficiently large $N$ it follows from Littlewood's Conjecture simply. Here we mean the
Littlewood's Conjecture proved by Konyagin \cite{Kon} and independently by McGehee, Pigno,
and Smith \cite{McGPigSmi} in 1981. See also pages 285-288 in \cite{DeVLor} for a book proof.
It is not difficult to prove the statement in general even in the case $n_1 = 0$ without
using Littlewood's Conjecture. One possible way is to use the identity
$$\sum_{j=1}^{n_N}{T\left(\frac{(2j-1)\pi}{n_N}\right)} = 0\,.$$
See \cite{KonLev}, for example.
Another way is to use Theorem 2 of \cite{Mer}. So there is certainly no shortage of possible
approaches to prove the starting observation of this paper even in the case $n_1 = 0$.

It seems likely that the number of zeros of the above sums in a period
must tend to $\infty$ with $N$. In a private communication Conrey asked
how fast the number of real zeros of the above sums in a period tends to $\infty$ as
a function $N$. In \cite{ConGranPooSou} the authors observed that
for an odd prime $p$ the Fekete polynomial
$$f_p(z)=\sum_{k=0}^{p-1}{\binom{k}{p} z^k}$$
(the coefficients are Legendre symbols) has $\sim \kappa_0 p$  zeros on the unit circle, where
$0.500813>\kappa_0>0.500668$. Conrey's question in general does not appear to be easy.

Littlewood in his 1968 monograph `Some Problems in Real and Complex Analysis´
\cite{Lit4}  poses the following research problem (problem~22),
which appears to still be open: `If the $n_m$ are integral and all
different, what is the lower bound on the number of real zeros of
$\sum_{m=1}^N \cos (n_m\theta)$? Possibly $N-1$, or not much less.´
Here real zeros are counted in a period.
In fact no progress appears to have been made on this
in the last half century. In a recent paper \cite{BorErdFerLoc} we showed that this is
false. There exist cosine polynomials $\sum_{m=1}^N \cos (n_m\theta)$ with
the $n_m$ integral and all different so that the number of its real zeros
in a period is $O(N^{9/10}(\log N)^{1/5})$ (here the frequencies $n_m = n_m(N)$
may vary with $N$). However, there
are reasons to believe that a cosine polynomial $\sum_{m=1}^N \cos (n_m\theta)$
always has many zeros in a period."

Let
$${\mathcal L}_n := \left\{P: \, P(z)=\sum_{j=0}^n{a_jz^j}\,, \enskip a_j \in \{-1,1\} \right\} \,.$$
Elements of ${\mathcal L}_n$ are often called Littlewood polynomials of degree $n$. Let
$${\mathcal K}_n := \left\{P: \, P(z)=\sum_{j=0}^n{a_jz^j}\,, \enskip a_j \in {\mathbb C}, 
\enskip |a_0|=|a_n|=1, \enskip |a_j| \leq 1 \right\}\,.$$
Observe that ${\mathcal L}_n \subset {\mathcal K}_n$.
In \cite{BorErdLit} we proved that any polynomial $P \in {\mathcal K}_n$
has at least $8n^{1/2}\log n$ zeros in any open disk centered at a point on the unit circle
with radius $33n^{-1/2}\log n$. Thus polynomials in ${\mathcal K}_n$
have quite a few zeros near the unit circle. One may naturally ask how many unimodular roots a
polynomial in ${\mathcal K}_n$ can have.
Mercer \cite{Mer} proved that if a Littlewood polynomial $P \in {\mathcal L}_n$ of the form (1.1)
is skew reciprocal, that is, $a_j = (-1)^ja_{n-j}$ for each $j=0,1,\ldots,n$, then it has
no zeros on the unit circle. However, by using different elementary methods it was observed
in both \cite{Erd1} and \cite{Mer} that if a Littlewood polynomial $P$ of the form (1.1) is self-reciprocal,
that is, $a_j = a_{n-j}$ for each $j=0,1,\ldots,n$, $n \geq 1$, then it has at least one zero
on the unit circle.
Mukunda \cite{Muk} improved this result by showing that every self-reciprocal Littlewood polynomial of
odd degree has at least $3$ zeros on the unit circle. Drungilas \cite{Dru} proved that
every self-reciprocal Littlewood polynomial of odd degree $n \geq 7$ has at least $5$ zeros
on the unit circle and every self-reciprocal Littlewood polynomial of even degree $n \geq 14$
has at least $4$ zeros on the unit circle. In [4] two types of Littlewood polynomials are considered:
Littlewood polynomials with one sign change in the sequence of coefficients and Littlewood polynomials
with one negative coefficient, and the numbers of the zeros such Littlewood polynomials have
on the unit circle and inside the unit disk, respectively, are investigated. Note that the
Littlewood polynomials studied in \cite {BorChoFerJan} are very special. In \cite{BorErdFerLoc} 
we proved that the average number of zeros of self-reciprocal Littlewood polynomials of degree 
$n$ is at least $n/4$. However, it is much harder to give decent lower bounds for the quantities
$$\text {\rm NZ}_n := \min_{P}{\text {\rm NZ}(P)}\,,$$
where $\text {\rm NZ}(P)$ denotes the number of zeros of a polynomial $P$ lying on the
unit circle and the minimum is taken for all self-reciprocal Littlewood polynomials
$P \in {\mathcal L}_n$. It has been conjectured for a long time that
$\lim_{n \rightarrow \infty}\text {\rm NZ}_n = \infty$.
In \cite{Erd8} we showed that
$\lim_{n \rightarrow \infty}\text {\rm NZ}(P_n) = \infty$
whenever $P_n \in {\mathcal L}_n$ is self-reciprocal and
$\lim_{n \rightarrow \infty}{|P_n(1)|} = \infty$.
This follows as a consequence of a more general result, see Corollary 2.3 in \cite{Erd8}, stated
as Corollary 1.5 here, in which the coefficients of the self-reciprocal polynomials
$P_n$ of degree at most $n$ belong to a fixed finite set of real numbers.
In \cite{BorErd3} we proved the following result.

\begin{thm}
If the set $\{a_j:j \in {\mathbb N}\} \subset {\mathbb R}$ is finite, the set
$\{j \in {\mathbb N}: a_j \neq 0\}$ is infinite, the sequence $(a_j)$ is not eventually periodic, and
$$T_n(t) = \sum_{j=0}^n{a_j \cos(jt)}\,,$$
then $\lim_{n \rightarrow \infty}{\text {\rm NZ}(T_n)} = \infty\,.$
\end{thm}

In \cite{BorErd3} Theorem 1.1 is stated without the assumption that the sequence
$(a_j)$ is not eventually periodic. However, as the following example shows,
Lemma 3.4 in \cite{BorErd3}, dealing with the case of eventually periodic sequences $(a_j)$,
is incorrect. Let
\begin{equation*}
\begin{split}
T_n(t) := & \cos t + \cos((4n+1)t) + \sum_{k=0}^{n-1}{(\cos((4k+1)t) - \cos((4k+3)t))} \\ 
= & \frac{1+\cos((4n+2)t)}{2\cos t} + \cos t\,. \\ 
\end{split}
\end{equation*}
It is easy to see that $T_n(t) \neq 0$ on $[-\pi,\pi] \setminus \{-\pi/2,\pi/2\}$ and the
zeros of $T_n$ at $-\pi/2$ and $\pi/2$ are simple. Hence $T_n$ has only two (simple) zeros
in a period. So the conclusion of Theorem 1.1 above is false for
the sequence $(a_j)$ with $a_0 := 0$, $a_1 := 2$, $a_3 :=-1$, $a_{2k} := 0$, $a_{4k+1} := 1$, $a_{4k+3} := -1$
for every $k=1,2,\ldots$. Nevertheless, Theorem 1.1 can be saved even
in the case of eventually periodic sequences $(a_j)$ if we
assume that $a_j \neq 0$ for all sufficiently large $j$. See Lemma 3.11 in \cite{Erd8} where
Theorem 1 in \cite{BorErd3} is corrected as

\begin{thm}
If the set $\{a_j:j \in {\mathbb N}\} \subset {\mathbb R}$ is finite,
$a_j \neq 0$ for all sufficiently large $j$, and
$$T_n(t) = \sum_{j=0}^n{a_j \cos(jt)}\,,$$
then $\lim_{n \rightarrow \infty}{\text {\rm NZ}(T_n)} = \infty\,.$
\end{thm}

It was expected that the conclusion of the above theorem remains true
even if the coefficients of $T_n$ do not come from the same sequence,
that is,
$$T_n(t) = \sum_{j=0}^n{a_{j,n} \cos(jt)}\,,$$
where the set
$$S := \{a_{j,n}: j \in \{0,1,\ldots,n\}, \, n \in {\mathbb C}\} \subset {\mathbb R}$$
is finite and
$$\lim_{n \rightarrow \infty}{|\{j \in \{0,1,\ldots,n\}, \, a_{j,n} \neq 0\}|} = \infty\,.$$

Associated with an algebraic polynomial
$$P(z) = \sum_{j=0}^n{a_jz^j}\,, \qquad a_j \in {\mathbb C}\,,$$
let
$$\text {\rm NC}_k(P) := \left|\left\{u: 0 \leq u \leq n-k+1, \enskip a_u + a_{u+1} + \cdots + a_{u+k-1} \neq 0 \right\}\right|\,.$$
In \cite{Erd8}  we proved the following results.

\begin{thm}
If $S \subset {\mathbb R}$ is a finite set, $P_{2n} \in {\mathcal P}_{2n}(S)$ are self-reciprocal polynomials,
$$T_n(t) := P_{2n}(e^{it})e^{-int}\,,$$
and
$$\lim_{n \rightarrow \infty}{\text {\rm NC}_k(P_{2n})} = \infty$$
for every $k \in {\mathbb N}$, then
$$\lim_{n \rightarrow \infty}{\text {\rm NZ}(P_{2n})} = \lim_{n \rightarrow \infty}{\text {\rm NZ}(T_n)} = \infty\,.$$
\end{thm}

Some of the most important consequences of the above theorem obtained in \cite{Erd8} are stated below.

\begin{cor}
If $S \subset {\mathbb R}$ is a finite set, $P_n \in {\mathcal P}_n(S)$ are self-reciprocal polynomials,
and
$$\lim_{n \rightarrow \infty}{|P_n(1)|} = \infty\,,$$
then
$$\lim_{n \rightarrow \infty}{\text {\rm NZ}(P_n)} = \infty\,.$$
\end{cor}

\begin{cor}
Suppose the finite set $S \subset {\mathbb R}$ has the property that
$$s_1 + s_2 + \cdots + s_k = 0\,, \enskip s_1,s_2,\ldots,s_k \in S\,, 
\enskip \text {\rm implies} \enskip s_1=s_2= \cdots =s_k=0\,,$$
that is, any sum of nonzero elements of $S$ is different from $0$.
If $P_n \in {\mathcal P}_n(S)$ are self-reciprocal polynomials and
$$\lim_{n \rightarrow \infty}{\text {\rm NC}(P_n)} = \infty\,,$$
then
$$\lim_{n \rightarrow \infty}{\text {\rm NZ}(P_n)} = \infty\,.$$
\end{cor}

J. Sahasrabudhe \cite{Sah} examined the case when $S \subset {\mathbb Z}$ is finite.
Exploiting the assumption that the coefficients are integer he proved
that for any finite set $S \subset {\mathbb Z}$ a self-reciprocal polynomial
$P \in {\mathcal P}_{2n}(S)$ has at least
$$c \, \left( \log\log\log|P(1)| \right)^{1/2-\varepsilon}-1$$
zeros on the unit circle of ${\mathbb C}$ with a constant $c > 0$ depending only on
$M = M(S) := \max\{|z|: z \in S\}$ and $\varepsilon > 0$.

Let $\phi(n)$ denote the Euler's totient function defined as the number of integers
$1 \leq k \leq n$ that are relative prime to $n$.
In an earlier version of his paper Sahasrabudhe \cite{Sah} used the trivial estimate
$\phi(n) \neq \sqrt n$ for $n \geq 3$ and he proved his result with the exponent
$1/4-\varepsilon$ rather than $1/2-\varepsilon$. Using the nontrivial estimate
$\phi(n) \geq n/(8\log\log n)$ in \cite{RosSch} for all $n > 3$ allowed him to prove
his result with $1/2-\varepsilon$.

In the papers \cite{BorErd3}, \cite{Erd8}, and \cite{Sah} the already mentioned Littlewood Conjecture, 
proved by Konyagin \cite{Kon} and independently by McGehee, Pigno, and Smith \cite{McGPigSmi}, plays a 
key role, and we rely on it heavily in the proof of the main results of this paper as well.
This states the following.

\begin{thm}
There is an absolute constant $c > 0$ such that
$$\int_{0}^{2\pi}{\Big| \sum_{j=1}^m{a_je^{i\lambda_jt}} \Big| \, dt} \geq c\gamma \log m\,$$
whenever $\lambda_1, \lambda_2, \ldots, \lambda_n$ are distinct integers and
$a_1, a_2, \ldots, a_m$ are complex numbers of modulus at least $\gamma > 0$.
Here $c=1/30$ is a suitable choice.
\end{thm}

This is an obvious consequence of the following result a book proof of which has been
worked out by Lorentz and DeVore, see pages 285--288 in \cite{DeVLor}.

\begin{thm}
If $\lambda_1 < \lambda_2 < \cdots < \lambda_m$ are integers and $a_1, a_2, \ldots, a_m$
are complex numbers, then
$$\int_{0}^{2\pi}{\Big| \sum_{j=1}^m{a_je^{i\lambda_j t}} \Big| \, dt} \geq
\frac{1}{30} \sum_{j=1}^m{\frac{|a_j|}{j}}\,.$$
\end{thm}

Associated with a finite set $S \subset {\mathbb C}$ we will use the notation $M(S) := \max\{|z|: z \in S\}$
throughout the paper.

\medskip

\section{New Results.}

The goal of this paper is to improve the exponent $1/2-\varepsilon$ to $1 -\varepsilon$ in Sahasrabudhe's
lower bound in \cite{Sah}  mentioned in Section 1. Sahasrabudhe's new idea is combined with the approach used in 
\cite{Erd8} offering an essencially simplified way to achieve our improvement.

Let, as before, $\text{\rm NZ}(T)$ denote the number of real zeros (by counting multiplicities) of a real trigonometric
polynomial $T$ in a period (of length $2\pi$).
Let $\text{\rm NZ}^*(T)$ denote the number of sign changes of a real trigonometric polynomial $T$ in a period 
(of length $2\pi$). Obviously $\text{\rm NZ}(T) \geq \text{\rm NZ}^*(T)$.

\begin{thm}
If $S \subset {\mathbb Z}$ is a finite set, $M = M(S) := \max\{|z|: z \in S\}$, $P \in {\mathcal P}_{2n}(S)$ is a self-reciprocal
polynomial,
$$T(t) := P(e^{it})e^{-int}\,,$$
then
$$\text {\rm NZ}^*(T_n) \geq \left( \frac{c}{1 + \log M} \right) \frac{\log\log\log|P(1)|}{\log\log\log\log|P(1)|}-1$$
with an absolute constant $c > 0$, whenever the right-hand side is defined.
\end{thm}

Let, as before, $\text{\rm NZ}(P)$ denote the number of real zeros (by counting multiplicities) of an algebraic
polynomial $P$ on the unit circle.

\begin{cor}
If $S \subset {\mathbb Z}$ is a finite set, $M = M(S) := \max\{|z|: z \in S\}$, $P \in {\mathcal P}_n(S)$ is a self-reciprocal
polynomial, then
$$\text {\rm NZ}(P) \geq \left( \frac{c}{1 + \log M} \right) \frac{\log\log\log|P(1)|}{\log\log\log\log|P(1)|}-1$$
with an absolute constant $c > 0$, whenever the right-hand side is defined.
\end{cor}

This improves the exponent $1/2 - \varepsilon$ to $1 - \varepsilon$ in a recent breakthrough result \cite{Sah} by
Julian Sahasrabudhe. We note that in both Sahasrabudhe's paper and this paper the assumption that the finite set
$S$ contains only integers is deeply exploited. Our next result is an obvious consequence of Corollary 2.2.

\begin{cor}
If the set $S \subset {\mathbb Z}$ is finite, $M = M(S) := \max\{|z|: z \in S\}$,
$$T(t) = \sum_{j=0}^n{a_j\cos(jt)}\,, \qquad a_j \in S\,,$$
then
$$\text {\rm NZ}^*(T) \geq \left( \frac{c}{1 + \log M} \right) \frac{\log\log\log|T(0)|}{\log\log\log\log|T(0)|}-1$$
with an absolute constant $c > 0$, whenever the right-hand side is defined.
\end{cor}

\medskip

\section{Lemmas.}

Our first four lemmas are quite similar to some of the lemmas used in \cite{Erd8}, but some modifications in the
formulation of these lemmas and their proofs are needed.

\begin{lem}
If $S \subset {\mathbb C}$ is a finite set, $k \in {\mathbb N}$,
$$S_k := \{s_1 + s_2 + \cdots + s_k: \, s_j \in S \cup \{0\}\}\,,$$
$$\gamma := \min_{z \in S_k \setminus \{0\}}{|z|}\,, \qquad M := M(S) := \max\{|z|: z \in S\}\,,$$
$P \in {\mathcal P}_{2n}(S)$, $H(z) := z^k - 1$,
and
\begin{equation}
\text {\rm NC}(PH) \leq \mu\,, 
\end{equation}
then
$$\int_{-\delta}^{\delta}{|P(e^{it})| \, dt} > 
\frac{\gamma}{30k} \, \log(\text {\rm NC}_k(P)) - \frac{\pi^2\mu M}{\delta}$$
for every $\delta \in (0,\pi)$.
\end{lem}

\textbf{Proof of Lemma 3.1.}\\
We define
$$G(z) := \sum_{j=0}^{k-1}{z^j}$$ so that $H(z) = G(z)(z-1)$.
As $P \in {\mathcal P}_{2n}(S)$ and the set $S$ is finite, the set $S_k$ is also finite.
By Theorem 1.6 there is an absolute constant $c > 0$ such that
\begin{equation}
\int_0^{2\pi}{|(PG)(e^{it})| \, dt} \geq \frac{\gamma}{30} \, \log(\text {\rm NC}(PG)) 
\geq \frac{\gamma}{30} \, \log(\text {\rm NC}_k(P))\,. 
\end{equation}
We define
$$M_k := M(S_k) = \max\{|z|: z \in S_k\} \leq k\max\{|z|: z \in S\} \leq kM\,.$$
Combining this with (3.1) we have
\begin{equation*}
\begin{split}
|(PG)(e^{it})| = &\frac{1}{|e^{it}-1|} \, |(PH)(e^{it})| \leq \frac{\mu M_k}{|e^{it}-1|} \\ 
= & \frac{\mu kM}{|2\sin(t/2)|} \leq \frac{\pi\mu kM}{|2t|}\,, \qquad t \in (-\pi,\pi)\,. \\
\end{split}
\end{equation*}
It follows that
\begin{equation}
\int_{[-\pi,\pi] \setminus [-\delta,\delta]}{|(PG)(e^{it})| \, dt}  
\leq 2\pi \frac{\pi\mu kM}{2\delta} = \frac{\pi^2\mu kM}{\delta} \,.
\end{equation}
Now (3.2) and (3.3) give
\begin{equation*}
\begin{split} \int_{-\delta}^{\delta}{|P(e^{it})| \, dt} & \geq 
\frac 1k \int_{-\delta}^{\delta}{|(PG)(e^{it})| \, dt} \\ 
& = \frac 1k \left( \int_{0}^{2\pi}{|(PG(e^{it})| \, dt} - 
\int_{[-\pi,\pi] \setminus [-\delta,\delta]}{|(PG)(e^{it})| \, dt} \right) \\ 
& \geq \frac{\gamma}{30k} \, \log(\text {\rm NC}_k(P))  - \frac{\pi^2\mu kM}{k\delta} \\ 
& \geq \frac{\gamma}{30k} \, \log(\text {\rm NC}_k(P))  - \frac{\pi^2\mu M}{\delta}\,. 
\end{split}
\end{equation*}
\qed

\begin{lem}
If $S \subset {\mathbb R}$ is a finite set, $P \in {\mathcal P}_{2n}(S)$ is self-reciprocal, $k \in {\mathbb N}$,
$H(z) := z^k - 1$, (3.1) holds,   
$$T(t) := P(e^{it})e^{-int}\,, \qquad R(x) := \int_0^x{T(t) \, dt}\,,$$
$M := M(S) = \max\{|z|:z \in S\}$,  and $0 < \delta \leq (2k)^{-1}$, then
$$\max_{x \in [-\delta,\delta]}{|R(x)|} < 42k(\mu+1)M\,.$$
\end{lem}

\textbf{Proof of Lemma 3.2.}\\
Let
$$P(z) = \sum_{j=0}^{2n}{a_j z^j}\,, \qquad a_j \in S\,,$$
be self-reciprocal. We have  
$$T(t) = a_n + \sum_{j=1}^n{2a_{j+n}\cos(jt)}\,, \qquad a_j \in S\,.$$
Observe that (3.1) implies that
\begin{equation}
|\{j: n+k \leq j \leq 2n, \, a_j \neq a_{j-k}\}| \leq \mu \,.
\end{equation}
We have
$$R(x) = a_nx + \sum_{j=1}^n{\frac{2a_{j+n}\sin(jx)}{j}}\,.$$
Now (3.4) implies that
$$R(x) = a_0x + \sum_{m=1}^{u}{F_{m,k}(x)}\,,$$
where
$$F_{m,k}(x) := \sum_{j=0}^{n_m-1}{\frac{2A_{m,k} \sin((j_m+jk)x)}{j_m + jk}}$$
with some $A_{m,k} \in S$, $m=1,2,\ldots,u$, $j_m \in {\mathbb N}$, and $n_m \in {\mathbb N}$, where $u \leq k(\mu+1)$
(we do not know much about $j_m$ and $n_m$). Since $S \subset [-M,M]$, it is sufficient to prove that
$$\max_{x \in [-\delta,\delta]}{|F_{m,k}(x)|} \leq \frac{41}{2}\,, \qquad m=1,2,\ldots,u\,,$$
that is, it is sufficient to prove that if $j_0 \in {\mathbb N}$ and
$$F(x) := \sum_{j=0}^{\nu}{\frac{\sin((j_0+jk)x)}{j_0 + jk}}\,,$$
then
\begin{equation}
\max_{x \in [-\delta,\delta]}{|F(x)|} = \max_{x \in [0,\delta]}{|F(x)|} \leq \frac{41}{2}\,.
\end{equation} 
Note that the equality in (3.5) holds as $F$ is odd. To prove the inequality in (3.5) let
$x \in (0,\delta]$, where $0 < \delta \leq (2k)^{-1}$. We break the sum as
\begin{equation}
F = R + S\,, 
\end{equation}
where
$$R(x) := \sum_{\substack{j=0 \\ j_0 + jk \leq x^{-1}}}^{\nu}{\frac{\sin((j_0+jk)x)}{j_0 + jk}}$$
and
$$S(x) := \sum_{\substack{j=0 \\ x^{-1} < j_0 + jk}}^{\nu}{\frac{\sin((j_0+jk)x)}{j_0 + jk}}\,.$$
Here
\begin{equation} 
\begin{split} |R(x)| \leq \, & \sum_{\substack{j=0 \\ j_0 + jk \leq x^{-1}}}^{\nu}{\left| \frac{\sin((j_0+jk)x)}{j_0 + jk} \right|} 
\leq (x^{-1} + 1)|x| \leq 1 + |x| \\
\leq \, & 1 + \delta = 1 + (2k)^{-1} \leq \frac 32\,, \\ 
\end{split}
\end{equation}
where each term in the sum in the middle is estimated by
$$\left| \frac{\sin((j_0+jk)x)}{j_0 + jk} \right| \leq \left| \frac{(j_0+jk)x)}{j_0 + jk} \right| = |x|\,,$$
and the number of terms in the sum in the middle is clearly at most $x^{-1} + 1$.
Further, using Abel rearrangement, we have
$$S(x) = -\frac{B_v(x)}{j_0 + vk} + \frac{B_u(x)}{j_0 + uk} + \sum_{\substack{j=0 \\ x^{-1} < j_0 + jk}}^{\nu}
{B_j(x)}\left(\frac{1}{j_0+jk} - \frac{1}{j_0+(j+1)k} \right)$$
with
$$B_j(x) := B_{j,k}(x) := \sum_{h = 0}^j{\sin((j_0 + hk)x)}$$
and with some $u,v \in {\mathbb N}_0$ for which $x^{-1} < j_0 + (u+1)k$ and $x^{-1} < j_0 + (v+1)k$.
Hence,
\begin{equation}
\begin{split}
|S(x)| \leq & \left| \frac{B_v(x)}{j_0 + vk} \right| + \left| \frac{B_u(x)}{j_0 + uk} \right| \\ 
+ & \sum_{\substack{j=0 \\ x^{-1} < j_0 + jk}}^{\nu}{|B_j(x)}|\left(\frac{1}{j_0+jk} - \frac{1}{j_0+(j+1)k} \right)\,. \\
\end{split}
\end{equation}
Note that $x \in (0,\delta]$, $0 < \delta \leq (2k)^{-1}$, $x^{-1} < j_0 + (w+1)k$, and $w \in {\mathbb N}_0$ imply
$$x^{-1} < j_0 + (w+1)k < 2(j_0 + wk) \qquad \text {\rm if} \enskip w \geq 1\,,$$
and
$$2k \leq \delta^{-1} \leq x^{-1} < j_0 + k \qquad \text {\rm if} \enskip w = 0\,,$$
and hence
\begin{equation}
\frac{1}{j_0 + wk} \leq 2x\,.
\end{equation} 
Observe also that $x \in (0,\delta]$ and $0 < \delta \leq (2k)^{-1}$ imply that $0 < x < \pi k^{-1}$.
Hence, with $z = e^{ix}$ we have
\begin{equation}
\begin{split} |B_j(x)| & = \left| \text {\rm Im} \left(\sum_{h=0}^j{z^{j_0 + hk}} \right) \right|
\leq  \left| \sum_{h=0}^j{z^{j_0 + hk}}\right| = \left| \sum_{h=0}^j{z^{hk}}\right| \\ 
& = \left| \, \frac{1-z^{(j+1)k}}{1-z^k} \right| \leq |1-z^{(j+1)k}| \, \frac{1}{|1-z^k|} \leq \frac{2}{|1-z^k|} \\ 
& \leq \frac{2}{2\sin(kx/2)} \leq \frac{\pi}{kx}\,. \\
\end{split}
\end{equation}
Combining (3.8), (3.9), and (3.10), we conclude
\begin{equation}
|S(x)| \leq \frac{\pi}{kx} \, 2x +  \frac{\pi}{kx} \, 2x +  \frac{\pi}{kx} \, 2x \leq \frac{6\pi}{k}\,. 
\end{equation}
Now (3.6), (3.7), and (3.11) give the inequality in (3.5) as $3/2 + 6\pi/k \leq 41/2$.
\qed

Our next lemma was used in \cite{Erd8}  in the same form. To prove it by contradiction is a simple exercise.
        
\begin{lem}
If $R$ is a continuously differentiable real-valued function on the interval $[-\delta,\delta]$, $\delta > 0$,
$$L := \int_{-\delta}^{\delta}{|R^{\prime}(x)| \, dx}\qquad and \qquad N := \max_{x \in [-\delta,\delta]}{|R(x)|}\,,$$
then there is an $\eta \in [-N,N]$ such that $R-\eta$ has at least $L(2N)^{-1}$ distinct zeros in $[-\delta,\delta]$.
\end{lem}

\begin{lem}
If $S \subset {\mathbb R}$ is a finite set, $k \in {\mathbb N}$, 
$$S_k := \{s_1 + s_2 + \cdots + s_k: \, s_j \in S \cup \{0\}\}\,,$$
$$\gamma := \min_{z \in S_k \setminus \{0\}}{|z|}\,, \qquad M := M(S) := \max\{|z|: z \in S\}\,,$$
$P \in {\mathcal P}_{2n}(S)$ is self-reciprocal, $T(t) := P(e^{it})e^{-int}$, $H(z) := z^k - 1$, and  
(3.1) holds, that is, 
$$\text {\rm NC}(PH) \leq \mu\,,$$ 
then 
$$\text {\rm NZ}^*(T) \geq
\left( \frac{\gamma}{30} \, \log(\text {\rm NC}_k(P)) - 2k\pi^2 \mu M \right)(84k(\mu+1)M)^{-1}\,.$$
\end{lem}

\textbf{Proof of Lemma 3.4.}\\
Let $0 < \delta := (2k)^{-1}$. Let $R$ be defined by
$$R(x) := \int_0^x{T(t) \, dt}\,.$$
Observe that $|T(x)| = |P(e^{ix})|$ for all $x \in {\mathbb R}$, and hence Lemma 3.1 yields that
\begin{equation*}
\begin{split} \int_{-\delta}^{\delta}{|R^\prime(x)| \, dx} & \, = \int_{-\delta}^{\delta}{|T(x)| \, dx}
= \, \int_{-\delta}^{\delta}{|P(e^{ix})| \, dx} \cr
& \, > \frac{\gamma}{30k} \, \log(\text {\rm NC}_k(P)) - \frac{\pi^2\mu M}{\delta} \cr
& \, = \frac{\gamma}{30k} \, \log(\text {\rm NC}_k(P)) - 2k\pi^2\mu M\,, \cr
\end{split}
\end{equation*}
while by Lemma 3.2 we have
$$\max_{x \in [-\delta,\delta]}{|R(x)|} < 42k(\mu+1)M \,.$$  
Therefore, by Lemma 3.3 there is an $\eta \in {\mathbb R}$ such that
$R - \eta$ has at least 
$$\left( \frac{\gamma}{30k} \, \log(\text {\rm NC}_k(P)) - 2k\pi^2 \mu M \right)(84k(\mu+1)M)^{-1}\,.$$
distinct zeros in $[-\delta,\delta]$. However, $T(x) = (R - \eta)^{\prime}(x)$ for all $x \in {\mathbb R}$, 
and hence 
$$\text {\rm NZ}^*(T) \geq \left( \frac{\gamma}{30} \, \log(\text {\rm NC}_k(P)) - 2k\pi^2\mu M  \right)(84k(\mu+1)M)^{-1}$$
follows by Rolle's Theorem.
\qed 

The following lemma, in which the assumption $S \subset {\mathbb Z}$ is crucial, is simple  
to prove. It is stated as Lemma 9 in \cite{Sah}. Its straightforward proof given in \cite{Sah} is reduced 
to the fact that a determinant of integer entries is an integer, and hence if it is not $0$, 
then its modulus is at least $1$.  

\begin{lem}
For $b \in {\mathbb N}$ let $A$ be a $b \times b$ invertible matrix with entries from $\widetilde{S} \subset {\mathbb Z}$.
If $A\textbf{x} = \textbf{y}$ with
$$\textbf{x} = (x(1),x(2),\ldots,x(d)) \in {\mathbb C}^b \quad \text{\rm and} \quad 
\textbf{y} = (y(1),y(2),\ldots,y(d)) \in {\mathbb C}^b\,,$$
then
$$\max\{|x(1)|,|x(2)|,\ldots,|x(d)|\} \leq M^{d-1}d^{d/2} \max\{|y(1)|,|y(2)|,\ldots,|y(b)|\}\,,$$
where $M := M(\widetilde{S}) := \max\{|z|: z \in \widetilde{S}\}$.
\end{lem}

For integers $1 \leq b \leq N$ we call
$$(x(1+r),x(2+r),\ldots,x(b+r))\in {\mathbb C}^N\,, \qquad r=0,1,\ldots,N-b\,,$$
the $b$-tuples of
$$(x(1),x(2),\ldots,x(N))\in {\mathbb C}^N\,.$$

The following lemma is Lemma 10 in \cite{Sah}.

\begin{lem}
For $u,v,b,t \in {\mathbb N}$,  let $S \subset {\mathbb Z}$ be a finite set such that $v-u > |S|^b+3b$, and let
$$(x(u+1),x(u+2),\ldots,x(v))\in S^{v-u}\,.$$
Let $V$ denote the linear space spanned by the $b$-tuples
$$(x(r+1),x(r+2),\ldots,x(r+b)) \in S^b\,, \qquad r=u,u+1,\ldots,v-b\,,$$
over ${\mathbb R}$. If $\text{\rm dim}(V) = t < b$, then there are
$$(x_j(u+b),x_j(u+b+1),\ldots,x_j(v-b)) \in {\mathbb C}^{v-u-2b+1}\,, \qquad j=1,2,\ldots, t\,,$$
such that
$$x(r) = x_1(r) + x_2(r) + \cdots + x_t(r)\,, \qquad r \in [u+b,v-b]\,,$$
where
$$(x_j(u+b),x_j(u+b+1),\ldots,x_j(v-b)) \in {\mathbb C}^{v-u-2b+1}$$
are periodic with period $\alpha_j \leq 16t\log\log(t+3)$ for each $j=1,2,\ldots,t$.
\end{lem}
Let $f$ be a continuous, even, real-valued function on $K := {\mathbb R}$ (mod $2\pi$) which
changes sign on $(0,\pi)$ exactly at $t_1 < t_2 < \cdots < t_d$, $d \geq 1$. We define the
{\it companion polynomial} $Q$ of $f$ by
$$Q(e^{it})e^{-id t} := (-1)^p2^d\prod_{j=1}^d{(\cos t - \cos t_j)}\,,$$
where $p \in \{0,1\}$ is chosen so that $f(t)e^{-idt}Q(e^{it}) \geq 0$
for all $t \in (-\pi,\pi)$.  Observe that
$$Q(z) : = (-1)^p\prod_{j=1}^d{(z - e^{it_j})(z - e^{-it_j})} =: \sum_{j=0}^{2d}{b_jz^j}$$
is a monic self-reciprocal algebraic polynomial of degree $2d$ with real coefficients and with constant
term $1$. Observe that
\begin{equation}
\begin{split}
|b_j| = \left| \frac{1}{2\pi}\int_{0}^{2\pi}{Q(e^{it})e^{-ijt} \, dt} \right| 
\leq & \max_{t \in [0,2\pi}{|Q(e^{it})|} \leq 4^d\,, \\ 
& j=0,1,\ldots,2d\,. \\
\end{split}
\end{equation}
Associated with $m \in {\mathbb N}$ let
$d_m := \text{\rm LCM}(1,2,\ldots,m)$. It is shown in \cite{RosSch} that $d_m < 3^m$ holds for all 
$m \in {\mathbb N}$, and this upper bound will be useful for us later in this paper. We remark though that
$\lim_{m \rightarrow \infty}{d_m/e^m} = 1$ holds and it is equivalent to the Prime Number Theorem,
see \cite{Ten}.

\begin{lem} Suppose $S \subset {\mathbb Z}$ is a finite set,
$P \in {\mathcal P}_{2n}(S)$ is self-reciprocal, and
$T(t) := P(e^{it})e^{-int}$ has exactly $2d \geq 2$ sign changes in $(-\pi,\pi)$.
Let $Q$ be the companion polynomial of $T$ (so the degree of the monic self-reciprocal algebraic 
polynomial $Q$ is $2d$) and let
\begin{equation}
F(z) := P(z)(z^{d_m} - 1)^2Q(z) = \sum_{j=0}^{2n+2d_m+2d}{a_jz^j}\,,
\end{equation} 
where $m := \lfloor 32d\log\log(2d+3) \rfloor$, and let
$$\{j \in [0,2n+2d_m+2d]: a_j \neq 0\} = \{j_1 < j_2 < \cdots < j_q\}\,.$$
If $1 \leq r \leq s \leq q$ are integers and   
\begin{equation}
|a_{j_k}| < (4M)^{-2d}(2d+1)^{-d-1/2}\,, \qquad k \in [r,s]\,,
\end{equation}
holds with $M := M(S) := \max\{|z|: z \in S\}$, then
$$s-r < (|S|+2)^{4m+2} + 6d + 3\,,$$
where $|S|$ denotes the number of elements in the set $S$.
\end{lem}

\textbf{Proof of Lemma 3.7.}\\
Let
$$S^* := \{s_1-s_2: \, s_1,s_2 \in S \cup \{0\}\} \subset {\mathbb Z}$$
and
$$S^{**} := \{s_1-s_2: \, s_1,s_2 \in S^*\} \subset {\mathbb Z}\,.$$
Obviously
\begin{equation}
|S^*| \leq (|S|+1)^2, \qquad  |S^{**}| \leq |S^*|^2 \leq (|S|+1)^4\,, 
\end{equation}
and
\begin{equation}
M(S^{**}) \leq 2M(S^*) \leq 4M(S)\,. 
\end{equation}
Put $b := 2d+1$. Suppose to the contrary that
\begin{equation}
s-r \geq (|S|+2)^{4m+2} + 3b\,. 
\end{equation}
Let
$$G_1(z) := P(z)(z^{d_m}-1) := \sum_{j=0}^{2n+d_m}{\alpha_jz^j}\,, \qquad \alpha_j \in S^*\,,$$
and
\begin{equation}
G_2(z) := G_1(z)(z^{d_m}-1) = P(z)(z^{d_m}-1)^2 = \sum_{j=0}^{2n+2d_m}{\beta_jz^j}\,, 
\end{equation}
where
\begin{equation}
\beta_j := \alpha_{j-d_m}-\alpha_j \in S^{**}, \qquad j=0,1,\ldots,2n+2d_m\,, 
\end{equation}
with
$$\alpha_{j-d_m} := 0, \quad \alpha_{2n+d_m+j+1} := 0\,, \qquad j=0,1,\ldots,d_m-1\,.$$
Observe that $G_2 \in {\mathcal P}_{2n+2d_m}(S^{**})$.
Let $V$ denote the linear space spanned by the $b$-tuples
$$(\beta_j,\beta_{j+1},\ldots,\beta_{j+b-1})\,, \qquad j \in [j_r,j_s-b+1]\,,$$
over ${\mathbb R}$. Using Lemma 3.5 with $\widetilde{S} := S^{**}$, (3.14), (3.16), and the fact that the 
polynomial $Q$ of degree $2d$ is monic, we can deduce that 
$\text{\rm dim}(V) = t < b$. It follows from (3.15) and (3.17) that
$$j_s - j_r \geq s-r \geq (|S|+2)^{4m+2} + 3b > |S^{**}|^m + 3b\,.$$
As $t < b = 2d+1$, we have $\lfloor 16t\log\log(t+3) \rfloor \leq \lfloor 32d\log\log(2d+3) \rfloor = m$.
Applying Lemma 3.6 we obtain that
$$(\beta_{j_r+b},\beta_{j_r+b+1},\ldots,\beta_{j_s-b+1})$$
is periodic with period $d_m$, that is,
\begin{equation}
\beta_{r+d_m} = \beta_r\,, \qquad r \in [j_r+b,j_s-d_m-b+1]\,. 
\end{equation}
We claim that
\begin{equation}
\beta_u = 0, \qquad u \in [j_r+b,j_s-b-|S^*|d_m]\,. 
\end{equation}
Indeed, if $\beta_u \neq 0$ for some $u \in [j_r+b,j_s-b-|S^*|d_m]$, then (3.19) and (3.20) give
\begin{equation*}
\begin{split} \alpha_{u + hd_m} & = \alpha_u + \sum_{j=1}^h{(\alpha_{u+jd_m}-\alpha_{u+(j-1)d_m})} 
= \alpha_u - \sum_{j=1}^h{\beta_{u+jd_m}} \\ 
& = \alpha_u - h\beta_u \in S^*, \qquad h = 0,1,\ldots,|S^*|\,, \\ 
\end{split}
\end{equation*}
exhibiting $|S^*|+1$ distinct elements of $S^*$, which is impossible.
It follows from (3.13), (3.18), and (3.21) that
$$a_j = 0, \qquad j \in [j_r+2b,j_s-b-|S^*|d_m]\,,$$
hence, recalling $d_m < 3^m$ and (3.15), we obtain
\begin{equation*}
\begin{split} 
s-r  & \leq |S^*|d_m + 3b < |S^*| 3^m + 3b < (|S|+1)^2 (|S|+2)^{4m} + 3b \\  
& \leq (|S|+2)^{4m+2} + 3b\,, \\ 
\end{split}
\end{equation*}
which contradicts (3.17). In conclusion
$$s-r < (|S|+2)^{4m+2} + 3b = (|S|+2)^{4m+2} + 6d + 3\,.$$
\qed

\begin{lem} Under the assumptions of Lemma 3.7 we have
$$\log q \leq 60\pi(8M)^{2d+1}(2d+1)^{d+3/2}(|S|+2)^{4m+2} + 6d + 3)\,.$$
\end{lem}

\textbf{Proof of Lemma 3.8.}\\
Let $L := (|S|+2)^{4m+2} + 6d + 3$, and $r := \lfloor q/L \rfloor$.
Observe that
$$P(e^{it})e^{-int}(e^{id_mt}-1)^2 e^{-id_mt}Q(e^{it})e^{-idt}$$
is real and nonnegative for all $t \in {\mathbb R}$. Combining this with Theorem 1.7 and Lemma 3.7 we obtain
\begin{equation}
\begin{split}
& \, \int_0^{2\pi}{P(e^{it})e^{-int}(e^{id_mt}-1)^2e^{-id_mt}Q(e^{it})e^{-idt} \,dt} \\
= & \int_0^{2\pi}{|P(e^{it})(e^{id_mt}-1)^2Q(e^{it})| \,dt} \\
\geq & \frac{1}{30} \, \sum_{k=1}^q{\frac{|a_{j_k}|}{k}} 
\geq \frac{1}{30} \, \sum_{j=1}^r \sum_{k=(j-1)L+1}^{jL}{\frac{|a_{j_k}|}{k}} \\
\geq & \frac{1}{30}(4M)^{-2d}(2d+1)^{-d-1/2} \sum_{j=1}^r{\frac{1}{jL}} \\
\geq & \frac{1}{30}(4M)^{-2d}(2d+1)^{-d-1/2}L^{-1}\log(r+1)\,. \\  
\end{split} 
\end{equation}
On the other hand, using orthogonality, (3.18), $\beta_n+d_m-j \in S^{**}$, and (3.12) we have
\begin{equation}
\begin{split}
\int_0^{2\pi}{P(e^{it})e^{-int}(e^{id_mt}-1)^2 e^{-id_mt}Q(e^{it})e^{-idt} \,dt} = &
2\pi \sum_{j=-d}^d{\beta_{n+d_m-j}b_j} \\ 
\leq & 2\pi(2d+1)(4M)4^d\,. \\
\end{split}
\end{equation}
Combining (3.22) and (3.23) we conclude
$$\frac{1}{30}(4M)^{-2d}(2d+1)^{-d-1/2}L^{-1}\log(r+1) \leq 2\pi(2d+1)(4M)4^d\,,$$
and hence
\begin{equation*}
\begin{split} 
\log q & \leq \log(r+1) + \log L \leq 60\pi(8M)^{2d+1}(2d+1)^{d+3/2}L \\
& \leq 60\pi(8M)^{2d+1}(2d+1)^{d+3/2}((|S|+2)^{4m+2} + 6d + 3)\,. \\ 
\end{split}
\end{equation*}
\qed 

Our final lemma follows easily from Lemma 3.6.

\begin{lem}
If $S \subset {\mathbb Z}$ is a finite set,
$P \in {\mathcal P}_{2n}(S)$, $0 \not\equiv R$ is a polynomial of degree at most $u$ with real coefficients,  
$$\text {\rm NC}(PR) \leq q\,,$$
$v := \lfloor 16u\log\log(u+3) \rfloor$, $k := d_v = \text{\rm LCM}(1,2,\ldots,v)$, and $H(z) = z^k-1$, then
$$\text{\rm NC}(PH) \leq \mu := (q +1)(k+|S|^{u+1}+3(u+1)+2)\,.$$
\end{lem}

\medskip

\section{Proof of the New Results.}

\textbf{Proof of Theorem 2.1.}\\
Suppose $S \subset {\mathbb Z}$ is a finite set,
$P \in {\mathcal P}_{2n}(S)$ is self-reciprocal, and $T(t) := P(e^{it})e^{-int}$
has exactly $2d$ sign changes in $(-\pi,\pi)$.
Without loss of generality we may assume that $d \geq 2$ otherwise we study the self-reciprocal polynomial
$\widetilde{P} \in {\mathcal P}_{2n}(\widetilde{S}$ defined by $\widetilde{P}(z) := (z^2+1)P(z)$, where
$M(\widetilde{S}) = 2M(S)$ and $\widetilde{P}(1) := 2P(1)$.
Let $Q$ be the companion polynomial of $T$. Let
$$F(z) := P(z)(z^{d_m} - 1)^2Q(z) = \sum_{j=0}^{2n+2d_m+2d}{a_jz^j}\,,$$
where $m := \lfloor 32d\log\log(2d+3) \rfloor$. Let
$$\{j \in [0,2n+2d_m+2d]: a_j \neq 0\} =: \{j_1 < j_2 < \cdots < j_q\}\,.$$
Lemma 3.8 together with $|S|+2 \leq 2M(S)+3 = 2M+3 \leq 5M$ and $4m+2d+3 \leq 5m$ implies
$$\log q \leq 60\pi(8M)^{2d+1}(5M)^{4m+2}(2d+1)^{d+3/2} + 180\pi(8M)^{2d+1}(2d+1)^{d+5/2}\,,$$
and hence
\begin{equation}
\log q \leq 240\pi(8M)^{5m}(2d+1)^{d+5/2}\,. 
\end{equation}
Applying Lemma 3.9 with $u := 2d_m + 2d$, we have
\begin{equation}
\text{\rm NC}(PH) \leq \mu := (q+1)(d_v+|S|^{u+1}+3(u+1)+2)
\end{equation}
with $v := \lfloor 16u\log\log(u+3) \rfloor$, $k := d_v$, and $H(z) = z^k-1$. Observe that if
$S \subset {\mathbb Z}$ and
$$S_k := \{s_1 + s_2 + \cdots + s_k: \, s_j \in S \cup \{0\}\}\,,$$
then
$$1 \leq \gamma := \min\{|z|: z \in S_k \setminus \{0\}\}\,.$$
Lemma 3.4 gives
\begin{equation}
\log(\text {\rm NC}_k(P)) \leq 2520k^2(\mu+1)kM(2d) + 60k^2\pi^2 \mu M \,. 
\end{equation}
Using $m := \lfloor 32d\log\log(2d+3) \rfloor$, $u := 2d_m + 2d$,
$v := \lfloor 16u\log\log(u+3) \rfloor$, $k = d_v < 3^v$, $d_m < 3^m$, (4.1), (4.2), $|S| \leq 2M(S)+1 = 2M+1$, and the
inequality $a+b \leq ab$ valid for all $a \geq 1$ and $b \geq 1$, we obtain
\begin{equation*}
\begin{split}
& \log(2520k(\mu+1)M(2d) + 60k^2\pi^2 \mu M) \leq \log(2520k^2(\mu+1)M(2d+1)) \\
\leq & \, 8 + 2\log k + \log (\mu + 1) + \log M + \log(2d+1) \\
\leq & \, 8 + 2\log k + \log(q+1) + \log d_v + \log(|S|^{u+1}) + \log(3u+6) \\ 
& \enskip \, + \log M + \log(2d+1) \\
\leq & \, 8 + 2(\log 3)v + \log 2 + 240\pi(8M)^{8m}(2d+1)^{d+5/2} + (\log 3)v \\ 
& \enskip \, + (u+1)\log(2M+1) + \log(3u+6) + \log M + \log(2d+1) \\
\leq & \, 9 + 2(\log 3) 16u\log\log(u+3) + 240\pi(8M)^{5m}(2d+1)^{d+5/2} \\ 
& \enskip \, + \, (2 \cdot 3^m + 2d+1)\log(2M+1) \\ 
& \enskip \, + \, (\log 2 + 3(\log 3) m)\log(2d+1)\log(2M+1) + \log M + \log(2d+1)\,, \\ 
\end{split}
\end{equation*}
and hence
\begin{equation*}
\begin{split} 
& \log(2520k(\mu+1)M(2d) + 60k^2\pi^2 \mu M) \\
\leq & 9 + 3(\log 3) 16(2 \cdot 3^m + 2d)\log\log(2 \cdot 3^m + 2d) \\ 
+ & 240\pi(8M)^{5m}(2d+1)^{d+5/2} +  (2 \cdot 3^m + 2d+1)\log(2M+1) \\ 
+ & (\log 2 + (\log 3) m)\log(2d+1)\log(2M+1) + \log M + \log(2d+1)\,. \\ 
\end{split}
\end{equation*}
Combining this with (4.3) and $m := \lfloor 32d\log\log(2d+3) \rfloor$ gives that there is an absolute constant 
$c_1>0$ such that 
\begin{equation}
\log\log\log \text {\rm NC}_k(P) \leq c_1(d\log(d+1) + d\log\log(2d+1)(1+\log M))\,. 
\end{equation}
It is easy to see that
$$\text {\rm NC}_k(P) \geq \frac {k|P(1)| - k^2 M}{kM} = \frac {|P(1)|}{M} - k\,.$$
Therefore if $|P(1)| \geq 2kM$, then  $\text {\rm NC}_k(P) \geq \frac 12 |P(1)|$, and the theorem 
follows from (4.4) after a straightforward calculus. 
If $|P(1)| < 2kM$, then it follows from $k := d_v < 3^v$, $v := \lfloor 16u\log\log(u+3) \rfloor$, 
$u := 2d_m + 2d < 2 \cdot 3^m + 2d$, and $m := \lfloor 32d\log\log(2d+3) \rfloor$ that 
\begin{equation*}
\begin{split}
\log\log\log |P(1)| < & \log\log\log(2kM) \leq c_2(1+\log M)\log(d+1) \\ 
\leq & c_2(1+\log M)(2d+1)\,, \\
\end{split}
\end{equation*}
with an absolute constant $c_2 > 0$, and the theorem follows.
\qed 

\textbf{Proof of Theorem 2.2.}\\
Let $S \subset {\mathbb Z}$ be a finite set.
If $P \in {\mathcal P}_{2\nu}(S)$ is self-reciprocal, then the corollary follows from Theorem 2.1.
If $P \in {\mathcal P}_{2\nu+1}(S)$ is self-reciprocal, then
$\widetilde{P} \in {\mathcal P}_{2\nu+2}(S^*)$ defined by
$$\widetilde{P}(z) := (z+1)P(z) \in {\mathcal P}_{2\nu+2}(S^*)$$
is also self-reciprocal, where the fact that $S \subset {\mathbb Z}$ is finite implies that the set
$$S^* := \{s_1 + s_2: \, s_1, s_2 \in S \cup \{0\}\} \subset {\mathbb Z}\,,$$
is also finite. Observe also that
$$M(S^*) = \max\{|z|: z \in S^*\} = 2\max\{|z|: z \in S\} = 2M(S)$$
and
$$\widetilde{P}(1) = 2P(1)\,.$$
Hence applying Theorem 2.1 to
$\widetilde{P} \in {\mathcal P}_{2\nu+2}(S^*)$, we obtain the statement of the
corollary for $P \in {\mathcal P}_{2\nu+1}(S)$ from Theorem 2.1 again.
\qed

\textbf{Proof of Theorem 2.3.}\\
The corollary follows from Theorem 2.1 and the fact that for every
trigonometric polynomial $T$ of the form
$$T(t) := a_0 + \sum_{j=1}^n{a_j\cos(jt)}\,, \qquad a_j \in {\mathbb Z}\,,$$
there is a self-reciprocal algebraic polynomial $P$ of the form
$$P(z) = 2a_0 z^n + \, \sum_{j=1}^n{a_j (z^{n+j} + z^{n-j})} \,, \qquad a_j \in {\mathbb Z}\,,$$
such that
$$2T(t) := P(e^{it})e^{-int}\,.$$
\qed

\section{Acknowledgements.}
The author wishes to thank Stephen Choi for his reading earlier versions of this paper carefully,


\begin{thebibliography}{HD}




\normalsize
\baselineskip=17pt


\bibitem[1]{AndBla} V.V. Andrievskii and H-P. Blatt
\emph{Discrepancy of Signed Measures and Polynomial Approximation},
Springer, New York, 2002.

\bibitem[2]{BloPol} A. Bloch and G. P\'olya
\emph{On the roots of certain algebraic equations},
Proc. London Math. Soc. {33} (1932), 102--114.

\bibitem[3]{Bor} P. Borwein
\emph{Computational Excursions in Analysis and Number Theory},
Springer, New York, 2002.

\bibitem[4]{BorChoFerJan} P. Borwein, S. Choi, R. Ferguson, and J. Jankauskas 
\emph{On Littlewood polynomials with prescribed number of zeros inside the unit disk},
Canad. J. of Math. {67} (2015), 507--526.

\bibitem[5]{BorErd1} P. Borwein and T. Erd\'elyi
\emph{Polynomials and Polynomial Inequalities},
Springer, New York, 1995.

\bibitem[6]{BorErd2} P. Borwein and T. Erd\'elyi
\emph{On the zeros of polynomials with restricted coefficients},
Illinois J. Math. {41} (1997), 667--675.

\bibitem[7]{BorErd3} P. Borwein and T. Erd\'elyi
\emph{Lower bounds for the number of zeros of cosine polynomials: a problem of Littlewood},
Acta Arith. {128} (2007), 377--384.

\bibitem[8]{BorErdFerLoc} P. Borwein, T. Erd\'elyi, R. Ferguson, and R. Lockhart
\emph{On the zeros of cosine polynomials : solution to a problem of Littlewood},
Ann. of Math. (2) {167} (2008), 1109--1117.

\bibitem[9]{BorErdKos1} P. Borwein, T. Erd\'elyi, and G. K\'os
\emph{Littlewood-type problems on $[0,1]$},
Proc. London Math. Soc. [79] (1999), 22--46.

\bibitem[10]{BorErdKos2} P. Borwein, T. Erd\'elyi, and G. K\'os
\emph{The multiplicity of the zero at  $1$ of polynomials with constrained coefficients}, 
Acta Arith. {159} (2013) no. 4, 387--395.

\bibitem[11]{BorErdLit} P. Borwein, T. Erd\'elyi, and F. Littmann
\emph{Polynomials with coefficients from a finite set},
Trans. Amer. Math. Soc. {360} (2008), 5145--5154.

\bibitem[12]{Boyd} D.W. Boyd
\emph{On a problem of Byrnes concerning polynomials with restricted coefficients},
Math. Comp. {66} (1997), 1697--1703.

\bibitem[13]{BuhCleWolZal} H. Buhrman, R. Cleve, R. de Wolf, and C. Zalka
\emph{Bounds for small-error and zero-error quantum algorithms}
40th Annual Symposium on Foundations of Computer Science, New York, 1999, 
358--368, IEEE Computer Soc., Los Alamitos, CA.

\bibitem[14]{CasKal} P.G. Casazza and N.J. Kalton
\emph{Roots of complex polynomials and Weyl-Heisenberg frame sets},
Proc. Amer. Math. Soc. {130} (2002) no. 8, 2313--2318.

\bibitem[15]{ConGranPooSou} B. Conrey, A. Granville, B. Poonen, and K. Soundararajan
\emph{Zeros of Fekete polynomials}, Ann. Inst. Fourier (Grenoble) {50} (2000), 865--889.

\bibitem[16]{CooDut} J.M. Cooper and A.M. Dutle
\emph{Greedy Galois games},
Amer. Math. Monthly {120} (2013) no. 5, 441–451. 

\bibitem[17]{CopRiv} D. Coppersmith and T.J. Rivlin
\emph{The growth of polynomials bounded at equally spaced points}
SIAM J. Math. Anal. {23} (1992) no. 4, 970--983.

\bibitem[18]{CrootHart} E. Croot and D. Hart
\emph{h-fold sums from a set with few products},
SIAM J. Discrete Math. {24} (2010) no. 2, 505--519.

\bibitem[19]{DeVLor} R.A. DeVore and G.G. Lorentz
\emph{Constructive Approximation}, Springer-Verlag, Berlin, 1993.

\bibitem[20]{DoNguVu} Y. Do, H. Nguyen, and V. Vu
\emph{Real roots of random polynomials: expectation and repulsion}, 
Proc. London Math. Soc. {111} (2015) no. 6, 1231--1260.

\bibitem[21]{Dru} P. Drungilas
\emph{Unimodular roots of reciprocal Littlewood polynomials},  
J. Korean Math. Soc. {45} (2008) no. 3, 835--840.

\bibitem[22]{Dub1} A. Dubickas
\emph{On the order of vanishing at $1$ of a polynomial},
Lithuanian Math. J. {39} (1999), 365--370.

\bibitem[23]{Dub2} A. Dubickas
\emph{Three problems of polynomials of small measure},
Acta Arith. {98} (2001), 279--292.

\bibitem[24]{Dub3} A. Dubickas
{Polynomials with multiple roots at $1$}, 
Int. J. Number Theory {10} (2014) no. 2, 391--400.  

\bibitem[25]{DudSchul} M. Dudik and L.J. Schulman
\emph{Reconstruction from subsequences},
Journal of Combinatorial Theory, Series A {103} (2003) no. 2, 337--348.

\bibitem[26]{EdeKos} A. Edelman and E. Kostlan
\emph{How many zeros of a random polynomial are real?}, 
Bull. Amer. Math. Soc. (N.S.) {32} (1995), 1--37. Erratum: 
Bull. Amer. Math. Soc. (N.S.) {33} (1996), 325.

\bibitem[27]{Erd1} T. Erd\'elyi
\emph{On the zeros of polynomials with Littlewood-type coefficient constraints}, 
Michigan Math. J. {49} (2001), 97--111.

\bibitem[28]{Erd2} T. Erd\'elyi
\emph{Markov-Bernstein type inequalities for polynomials under Erd\H os-type constraints}, 
Paul Erd\H os and his Mathematics I, Bolyai Society Mathematical Studies, 11, 
G\'abor Hal\'asz, L\'aszl\'o Lov\'asz,  Dezs\H o Mikl\'os, and Vera T. S\'os (Eds.), 
Springer Verlag, New York, 2002, 219--239.

\bibitem[29]{Erd3} T. Erd\'elyi
\emph{Polynomials with Littlewood-type coefficient constraints},
Approximation Theory X: Abstract and Classical Analysis, Charles K. Chui,
Larry L. Schumaker, and Joachim St\"ockler (Eds.), Vanderbilt University Press,
Nashville, TN, 2002, 153--196.

\bibitem[30]{Erd4} T. Erd\'elyi
{Extensions of the Bloch-P\'olya theorem on the number of real zeros of polynomials}, 
Journal de th\'eorie des nombres de Bordeaux {20} (2008), 281--287.

\bibitem[31]{Erd5} T. Erd\'elyi
\emph{An improvement of the Erd\H os-Tur\'an theorem on the zero distribution of the zeros of polynomials},
C. R. Acad. Sci. Paris S\'er. I Math. {346} (2008), 267--270.  

\bibitem[32]{Erd6} T. Erd\'elyi
\emph{Pseudo-Boolean functions and the multiplicity of the zeros of polynomials},
Journal d'Analyse Math. {127} (2015), 91--108.

\bibitem[33]{Erd7} T. Erd\'elyi
\emph{Coppersmith-Rivlin type inequalities and the order of vanishing of polynomials at $1$},
Acta Arith. {172} (2016) no. 3, 271--284.

\bibitem[34]{Erd8} T. Erd\'elyi
\emph{The number of unimodular zeros of self-reciprocal polynomials with coefficients in a finite set}, 
Acta Arith. {176} (2016)  no. 2, 177--200.

\bibitem[35]{ErdOff} P. Erd\H os and A. C. Offord
\emph{On the number of real roots of a random algebraic equation},  
Proc. London Math. Soc. {6} (1956), 139--160.

\bibitem[36]{ErdTur} P. Erd{\H{o}}s and P. Tur\'an 
\emph{On the distribution of roots of polynomials}, 
Ann. Math. {57} (1950), 105--119. 

\bibitem[37]{Fer} Le Baron O. Ferguson
\emph{Approximation by Polynomials with Integral Coefficients},
Amer. Math. Soc., Rhode Island, 1980.

\bibitem[38]{FosKras} W. Foster and I. Krasikov
\emph{An improvement of a Borwein-Erd\'elyi-K\'os result},
Methods Appl. Anal. {7} (2000) no. 4, 605--614.

\bibitem[39]{Gunt} C.S. G\"unt\"urk
\emph{Approximation by power series with $\pm 1$ coefficients},
Int. Math. Res. Not. (2005) no. 26, 1601--1610.

\bibitem[40]{Han} D. Hanson
\emph{On the product of the primes}
Canad. Math. Bull. {15} (1972), 33--37.

\bibitem[41]{HarWright} G.H. Hardy and E.M. Wright
\emph{An Introduction to the Theory of Numbers},
Clarendon Press, Oxford, 1938.

\bibitem[42]{Hua} L.K. Hua
\emph{Introduction to Number Theory},
Springer-Verlag, Berlin Heidelberg, New York, 1982.

\bibitem[43]{Kac1} M. Kac
\emph{On the average number of real roots of a random algebraic equation}, 
Bull. Amer. Math. Soc. {49} (1943), 314--320.

\bibitem[44]{Kac2} M. Kac
\emph{On the average number of real roots of a random algebraic equation. II},  
Proc. London Math. Soc. {50} (1949), 390--408.

\bibitem[45]{Kon} S.V. Konyagin
\emph{On a problem of Littlewood}, 
Mathematics of the USSR, Izvestia {18} (1981), 205--225.

\bibitem[46]{KonLev} S.V. Konyagin and V.F. Lev
\emph{Character sums in complex half planes}, 
J. Theor. Nombres Bordeaux {16} (2004) no. 3, 587--606.

\bibitem[47]{KosLigSzik} G. K\'os, P. Ligeti, and P. Sziklai
\emph{Reconstruction of matrices from submatrices},
Math. Comp. {78} (2009), 1733--1747.

\bibitem[48]{Kras} I. Krasikov
\emph{Multiplicity of zeros and discrete orthogonal polynomials},
Results Math. {45} (2004), 1--2, 59--66.

\bibitem[49]{Lit1} J.E. Littlewood
\emph{On the mean values of certain trigonometrical polynomials},  
J. London Math. Soc. {36} (1961), 307--334.

\bibitem[50]{Lit2} J.E. Littlewood
\emph{On the real roots of real trigonometrical polynomials (II)},  
J. London Math. Soc. {39} (1964), 511--552.

\bibitem[51]{Lit3} J.E. Littlewood
\emph{On polynomials $\sum \pm z^m$ and $\sum e^{\alpha_{m}i} z^m$, $z=e^{\theta i}$},  
J. London Math. Soc. {41} (1966), 367--376.

\bibitem[52]{Lit4} J.E. Littlewood
\emph{Some Problems in Real and Complex Analysis},  
Heath Mathematical Monographs, Lexington, Massachusetts, 1968.

\bibitem[53]{LitOff1} J.E. Littlewood and A.C. Offord
\emph{On the number of real roots of a random algebraic equation. II},  
Proc. Cambridge Philos. Soc. {35} (1939), 133--148.

\bibitem[54]{LitOff2} J.E. Littlewood and A.C. Offord
\emph{On the number of real roots of a random algebraic equation. III},  
Rec. Math. [Mat. Sbornik] N.S. {54} (1943), 277--286,

\bibitem[55]{McGPigSmi} O.C. McGehee, L. Pigno, and B. Smith
\emph{Hardy's inequality and the $L_1$ norm of exponential sums},  
Ann. Math. {113} (1981), 613--618.

\bibitem[56]{Mer} I.D. Mercer
\emph{Unimodular roots of special Littlewood polynomials}, 
Canad. Math. Bull. {49} (2006) no. 3, 438--447.

\bibitem[57]{Moss} M.J. Mossinghoff
\emph{Polynomials with restricted coefficients and prescribed noncyclotomic factors},
London Math. Soc. J. Comput. Math. {6} (2003), 314--325 (electronic).

\bibitem[58]{Muk} K. Mukunda,
\emph{Littlewood Pisot numbers}, 
J. Number Theory {117} (2006) no. 1, 106--121.

\bibitem[59]{NguNguVu} H. Nguyen, O. Nguyen, and V. Vu
\emph{On the number of real roots of random polynomials},
Commun. Contemp. Math. {18} (2016), 1550052.

\bibitem[60]{OdPoo} A.M. Odlyzko and B. Poonen
\emph{Zeros of polynomials with $0,1$ coefficients},
Enseign. Math. (2) {39} (1993), 317--348.

\bibitem[61]{Pin} C. Pinner
\emph{Double roots of $[-1,1]$ power series and related matters},
Math. Comp. {68} (1999), 1149--1178.

\bibitem[62]{PriSol} I.E. Pritsker and A.A. Sola
\emph{Expected discrepancy for zeros of random algebraic polynomials},
Proc. Amer. Math. Soc. {142} (2014), 4251--4263.

\bibitem[63]{Rah} E.A. Rakhmanov 
\emph{Bounds for polynomials with a unit discrete norm},     
Ann. of Math. {165} (2007), 55--88.

\bibitem[64]{Rod} F. Rodier
\emph{Sur la non-lin\'earit\'e des fonctions bool\'eennes},  
Acta Arith. {115} (2004) no. 1, 1--22.

\bibitem[65]{RosSch} J.B. Rosser and L. Schoenfeld
\emph{Approximate formulas for some functions of prime numbers}, 
Illinois J. Math. {6} (1962), 64--94. 

\bibitem[66]{Sah} J. Sahasrabudhe
\emph{Counting zeros of cosine polynomials: on a problem of Littlewood}, 
Adv. Math., (2019), to appear. 

\bibitem[67]{Sch} E. Schmidt
\emph{\"Uber algebraische Gleichungen vom P\'olya-Bloch-Typos}, 
Sitz. Preuss. Akad. Wiss., Phys.-Math. Kl. (1932), 321.

\bibitem[68]{Schur} I. Schur
\emph{Untersuchungen \"uber algebraische Gleichungen},
Sitz. Preuss. Akad. Wiss., Phys.-Math. Kl. (1933), 403--428.

\bibitem[69]{Shpar} I.E. Shparlinski
\emph{Finite Fields: Theory and Computation: the Meeting Point of Number Theory, Computer Science,
Coding Theory and Cryptography},
Dordrecht, London, 1999.

\bibitem[70]{Sol} B. Solomyak
\emph{On the random series $\sum \pm \lambda^n$ (an Erd\H os problem)},
Ann. Math. {142} (1995), 611--625.

\bibitem[71]{Sze} G. Szeg\H o
\emph{Bemerkungen zu einem Satz von E. Schmidt uber algebraische Gleichungen},
Sitz. Preuss. Akad. Wiss., Phys.-Math. Kl. (1934), 86--98.

\bibitem[72]{TaoVu} T. Tao and V. Vu
\emph{Local universality of zeros of random polynomials},
Int. Math. Res. Notices {2015} (2015) no. 13, 5053--5139.

\bibitem[73]{Ten} G. Tenenbaum
\emph{Introduction to Analytic and Probabilistic Number Theory},
Cambridge University Press, Cambridge, England, 1995.

\bibitem[74]{TotVar} V. Totik and P. Varj\'u
\emph{Polynomials with prescribed zeros and small norm}, 
Acta Sci. Math. (Szeged) {73} (2007) no. 3-4, 593--611. 

\end{thebibliography}
\end{document}